\documentclass{article}
\usepackage{arxiv,xcolor}
\usepackage{multirow}
\usepackage{amsmath,amssymb}

\usepackage{subfigure, caption, mathtools}
\mathtoolsset{showonlyrefs}
\newcommand{\brac}[1]{\left(#1 \right)}
\newcommand{\pd}[2]{\frac{\partial{#1}}{\partial{#2} }}
\newcommand{\RR}{\mathbb R}

\title{Projection based model reduction for the immersed boundary method}

\author{
  Yushuang Luo \\
  Department of Mathematics \\
  The Pennsylvania State University, University Park, PA 16802, USA\\
  \texttt{yzl55@psu.edu} \\
  \And
  Xiantao Li \\
  Department of Mathematics, \\
    The Pennsylvania State University, University Park, PA 16802, USA\\
  \texttt{xxl12@psu.edu} \\
  \And
  Wenrui Hao \\
  Department of Mathematics\\
  The Pennsylvania State University, University Park, PA 16802, USA\\
  \texttt{wxh64@psu.edu} \\
}

\begin{document}
\maketitle



\begin{abstract}
Fluid-structure interactions are central to many bio-molecular processes, and they impose a great challenge for computational and modeling methods.  In this paper, we consider the immersed boundary method  (IBM) for biofluid systems, and to alleviate the computational cost,  we apply reduced-order techniques to eliminate the degrees of freedom associated with the large number of fluid variables. We show how reduced models can be derived using Petrov-Galerkin projection and subspaces that maintain the incompressibility condition.
More importantly, the reduced-order model is shown to  preserve the Lyapunov stability.  We also address the practical issue of computing coefficient matrices in the reduced-order model  using an interpolation technique.
The efficiency and robustness of the proposed formulation are examined with test examples from various applications. 
\end{abstract}


\section{Introduction}\label{sec:intro}

Biofluid dynamics, the study of cellular movement in biological fluid flow, is essential for  understanding how  the cellular behavior changes within living tissues \cite{kleinstreuer2006biofluid}. %
With the rapid development of scientific computing algorithms, mathematical modeling and numerical simulations have become an indispensable  approach for studying biofluid dynamics. Specifically, the interaction between cell structures and the surrounding fluid flow is of the {\color{black}utmost importance}. 
Mathematically speaking, this belongs to a large class of problems known as the fluid-structure interactions (FSI), often described  by coupling the incompressible Navier–Stokes equations with solid equations.
A variety of computational and modeling techniques  have been developed for FSIs, and they have been successfully implemented in
studying  (among many other applications)  biology and biomedical diseases \cite{Canic,CanM1,lei2012predicting,li2014probing,wang2002hydrodynamic,yang2009dynamic,YuYue,zhao2016decoupled}.
In order to numerically solve FSI problems, several numerical methods have been developed to represent/track the interface movement explicitly, such as the boundary element method (BEM) \cite{hall1994boundary,li2011boundary,everstine1990coupled}, the IBM \cite{peskin2002immersed,lai2000immersed,atzberger2006stochastic,sotiropoulos2014immersed}, the immersed interface method (IIM) \cite{leveque1994immersed}, the fictitious domain method (FDM) \cite{glowinski1994fictitious,hao2015fictitious}, and the front tracking method (FTM) \cite{tryggvason2001front,glimm2003conservative}. Another alternative approach to solve the FSI problem is  to capture the interface dynamics implicitly by evolving a scalar function defined on the whole domain. The level-set method \cite{cottet2006level}, the phase-field method \cite{du2004phase}, and the implicit boundary integral method \cite{kublik2013implicit} are important examples. 
However, direct simulations based on these methods tend to be time-consuming and computationally expensive for the prediction and analysis of long-term dynamics {\color{black}(although the short-term prediction is  certainly feasible)}. {\color{black}Often of interest in biology, is the structure dynamics, which, due to its observability, is easy to validate either experimentally or computationally \cite{VandV}. } 
In addition, there are many important scenarios where the cell structure is immersed in a large fluid environment, and simulating the entire system becomes computationally challenging.

The purpose of this paper is to explore an alternative to reduce the computational cost using reduced-order techniques,
which have been applied to a wide variety of problems in science and engineering \cite{Bai2002,benner2015survey,gugercin2008,anic2008interpolation}. 
Reduced-order modeling is concerned with large-dimensional dynamical systems with low-dimensional input and output, and the main objective is to construct reduced models that can approximate the mapping from the input directly to the output. {\color{black}The present FSI problem will be formulated as a reduced-order problem, where the input is the force exerted from the structure and the output is the local velocity of the structure.} {\color{black}As a proof-of-concept, we utilize the conventional IBM model \cite{peskin2002immersed}. Specifically, incompressible unsteady Stokes flows are considered, together with the no-slip {\color{black}interface} condition enforced on the immersed structure. But it is also important to point out that there have been many extensions of the original IBM framework with different treatments for the Lagrangian equations of motion or the fluid dynamics \cite{bao2017immersed,nitti2020immersed,GILMANOV2015814}, and reduced-order modeling can be considered in those settings as well. } Our starting point is a semi-discrete representation of the IBM model, so that {\color{black}the dynamics of fluid and structure motion can be expressed as coupled ODEs}, which can then be placed in the reduced-order modeling framework. Then we derive the effective mapping from the structure force to the local velocity, which completely eliminates the fluid variables. 

To construct specific reduced models that do not involve the fluid dynamics explicitly, we first construct subspaces that preserve the incompressibility condition, followed by a Petrov-Galerkin projection. We show that the choice of the subspaces ensures certain interpolation conditions on the underlying transfer function. 
An important departure from standard reduced-order problems is that in IBM, the structure is also evolving continuously. As a result, the subspaces are varying in time. {\color{black}This poses some challenges as the coefficient matrices of the reduced models need to be updated frequently.} To circumvent this issue, we observe the connection between the those matrices and the Green's function of the Laplace equation. More specifically, the entries of those matrices are tied to the nodal points on the structure. When the two points are far apart, the corresponding entry can be well approximated by the Green's function. On the other hand, for points that are within some cut-off distance, the computation can be done in advance, and then in the simulation, those entries can be computed by interpolation. We show that such a strategy avoids repetitive computation of those coefficient matrices and it can speed up the computation considerably.   

The remaining part of the paper is organized as follows: in Section \ref{sec:model}, we introduce the full-order model (FOM) in the IBM setup; in Section \ref{sec:reduce}, we formulate the reduced-order model (ROM); several numerical examples are used to compare both full-order and reduced-order models in Section \ref{sec:numres}; then the conclusion is drawn in Section \ref{sec:con}.

\section{Full-order Model}\label{sec:model}
In this section we briefly review the mathematical formulation of the IBM and derive its semi-discrete representation, which will serve as the full-order model (FOM).

\subsection{Mathematical formulation of the IBM}
The IBM is intended for the computer simulation of FSI, especially in biological fluid dynamics. It is mathematically defined by a set of differential equations involving a mixture of Eulerian and Lagrangian descriptions, linked by the Dirac delta function. {\color{black}The dynamics of the fluid is described in terms of the velocity $\boldsymbol{u}(\boldsymbol{x},t)$ and the pressure $p(\boldsymbol{x},t)$ on an Eulerian coordinate for $\boldsymbol{x}\in \Omega$, where $\Omega \subset \RR^d, \ d=2 \text{ or }3$, represents the fluid domain. The immersed structures, on the other hand, are handled in a Lagrangian coordinate as a parametric curve or surface $X(\boldsymbol{s},t)$. Specifically, $X(\boldsymbol{s},t)$ represents the position at time $t$ in Cartesian coordinates of the structure point labeled by $\boldsymbol{s}\in \Gamma$, where $\Gamma\subset \RR^{d-1}$ is the parameter space. In this work, we focus on two-dimensional flows where the structure is described as a parametric curve. In this case $\boldsymbol{x}\in\Omega \subset \RR^2$ and $s$ is a scalar parameter. The formulation is mostly algebraic. Therefore the extension to high dimensional cases is straightforward.} Assuming constant density, the time-dependent Stokes equation is used to model the incompressible flow
\begin{align}
  \label{eq:motion}
  \rho \pd{\boldsymbol{u}}{t} =  &- \nabla p + \mu \nabla^2 \boldsymbol{u} + \boldsymbol{f}, 
  \\
  \label{eq:mass}
  \nabla \cdot \boldsymbol{u} = & 0,
\end{align}
where $\rho$ and $\mu$ are the fluid density and viscosity, respectively. The body force $\boldsymbol{f}$ exerted by the structure on the fluid is defined as
\begin{equation}
  \label{eq:force}
  \boldsymbol{f}(\boldsymbol{x},t) = \int\limits_\Gamma \boldsymbol{F}(s,t) \, \delta(\boldsymbol{x} -
  \boldsymbol{X}(s,t)) \,ds, 
\end{equation}
where $\delta(\boldsymbol{x})$ is the Dirac delta function. In addition, $\boldsymbol{F}(s,t)$ denotes the force density on the immersed structure, defined as
\begin{gather}
  \label{eq:forceDensity}
  \boldsymbol{F}(s,t) = \boldsymbol{\mathcal{F}} \left[\boldsymbol{X}(s,t)\right],
\end{gather}
where $\boldsymbol{\mathcal{F}}$ is a functional of the IBM configuration. Spring forces, bending resistance or any other type of behavior (area and volume conservation constraints) can be built into this functional to embody the physics of the immersed structure under  different circumstances \cite{hao2015fictitious,pivkin2008accurate}. We give detailed description of the force density in section \ref{sec:numres} in the numerical examples for various systems. 

Assuming an over-damped structure, the immersed boundary must move with the local fluid velocity:
\begin{equation}
  \label{eq:membrane}
  \pd{\boldsymbol{X}(s,t)}{t} = \boldsymbol{u}(\boldsymbol{X}(s,t),t) = \int\limits_\Omega
  \boldsymbol{u}(\boldsymbol{x},t) \, \delta(\boldsymbol{x}-\boldsymbol{X}(s,t)) \, d\boldsymbol{x}.
\end{equation}
This last equation is nothing other than the no-slip condition written as a delta function convolution.

\subsection{The semi-discrete equations}
To derive a semi-discrete representation of the IBM, we use the finite difference discretization \cite{peskin2002immersed}. Other numerical methods  can also be applied to discretize the IBM, e.g., the finite element method \cite{boffi2003finite} and the finite volume method \cite{kim2001immersed}, in which the state space consists of nodal values.

\begin{figure}
    \centering
    \includegraphics[width=0.5\textwidth]{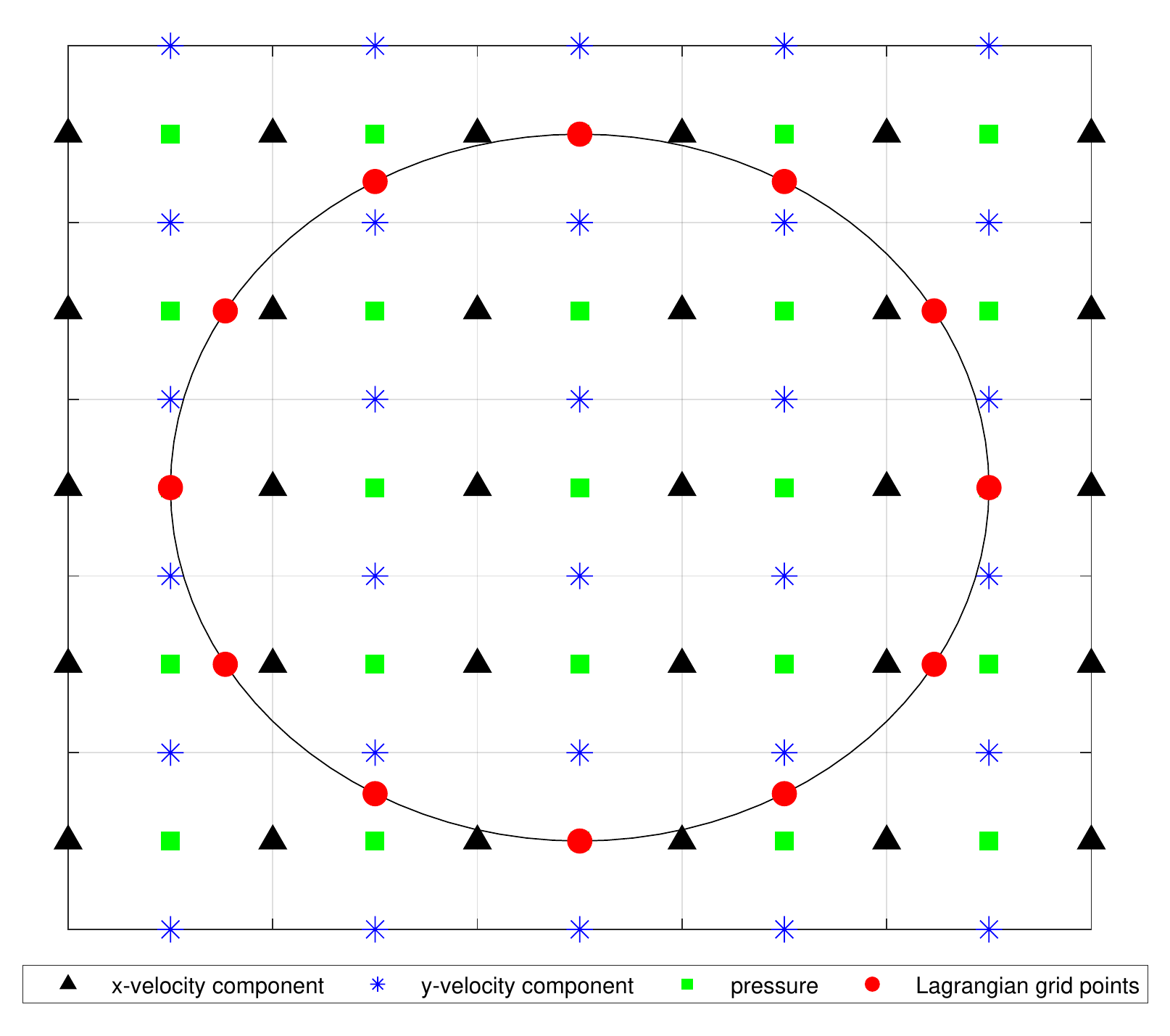}
    \caption{The fluid field variables are defined on a regular staggered grid. The structure variables are defined on an Lagrangian grid.}
    \label{fig: semidiscrete}
\end{figure}

In this work, fluid variables are discretized on a uniform staggered Eulerian grid, denoted $\Omega_h$; and the structures are discretized on an independent Lagrangian grid, denoted by $\Gamma_h$ (Fig \ref{fig: semidiscrete}). The Eulerian grid points are of the form $\boldsymbol{x} = \boldsymbol{j}h$, where $\boldsymbol{j} = (j_1,j_2)$ is a two-dimensional vector with integer components and $h$ is the Eulerian grid size. The Lagrangian grid is a set of $s$ of the form $k\Delta s$, where $k$ has integer components. The following restriction is imposed to avoid leak\cite{peskin2002immersed},
\begin{equation}\label{eq:avoidleak}
    \vert \boldsymbol{X}(s+\Delta s,t) - \boldsymbol{X}(s,t) \vert < \frac{h}{2},
\end{equation}
for all $s$.
First, the semi-discrete equations for \eqref{eq:motion}-\eqref{eq:mass} form a system of linear differential-algebraic equations (DAEs)
\begin{align}
    \rho\dot{\boldsymbol{u}}_h(t) = & -N\boldsymbol{p}_h(t) + \mu G\boldsymbol{u}_h(t) + \boldsymbol{f}_h(t), \label{eq:DAEmotion}\\
    0 = & M\boldsymbol{u}_h(t), \label{eq:DAEmass}
\end{align}
where 
\[
\boldsymbol{u}_h(t) = \begin{bmatrix}
\vdots\\
u_h(\boldsymbol{x},t)\\
\vdots
\end{bmatrix}\in \RR^{n_u}, \quad 
\boldsymbol{p}_h(t) = \begin{bmatrix}
\vdots\\
p_h(\boldsymbol{x},t)\\
\vdots
\end{bmatrix}\in \RR^{n_p}, \quad
\boldsymbol{f}_h(t) = \begin{bmatrix}
\vdots\\
f_h(\boldsymbol{x},t)\\
\vdots
\end{bmatrix}\in \RR^{n_u}, \quad \forall \boldsymbol{x}\in \Omega_h,
\]
are vectors of discrete velocity field, pressure and body force, respectively. $G \in \RR^{n_u,n_u}$ is the discrete Laplace operator. Matrices $N \in \RR^{n_u,n_p}$ and $M \in \RR^{n_p,n_u}$ are the discrete gradient and divergence operators, respectively. 

Second, the integrals in \eqref{eq:force} and \eqref{eq:membrane} are replaced by the following sums over the appropriate grid points,
\begin{align}
    f_h(\boldsymbol{x},t) = & \sum_{s\in\Gamma_h} F_h(s,t) \ \delta_r(\boldsymbol{x}-X_h(s,t)) \ \Delta s \quad \forall \boldsymbol{x}\in\Omega_h, \label{eq:sumf} \\
    \dot{X}_h(s,t) = & \sum_{\boldsymbol{x}\in\Omega_h} u_h(\boldsymbol{x},t) \ \delta_r(\boldsymbol{x}-X_h(s,t)) \ h^2 \quad \forall s\in\Gamma_h, \label{eq:sumv}
\end{align}
where $F_h(s,t)$ is the discrete Lagrangian force density associated with the structure point labeled $s$, obtained by discretizing \eqref{eq:forceDensity}. In addition, a function $\delta_r(\boldsymbol{x})$ that is nonsingular for each $r$ but approaches $\delta(\boldsymbol{x})$ as $r\to 0$ is needed. There are many ways to construct such $\delta_r$. We choose a radially symmetric function with compact support as follows \cite{yang2012generalized}, 
\begin{equation}\label{eq:deltar}
    \delta_r(\boldsymbol{x}) = \begin{cases}
    C_r\left[ 1+\frac{|\boldsymbol{x}|^2}{r^3}(2|\boldsymbol{x}|-3r)\right] & |\boldsymbol{x}| \leq r, \\
    0 & |\boldsymbol{x}|>r,
    \end{cases}
\end{equation}
where the normalizing constant $C_r = \frac{10}{3\pi}r^{-2}$ depends on $r$ and the space dimension ($C_r=\frac{15}{4\pi}r^{-3}$ in 3D).
For computational efficiency, we choose $r=2h$ in all our numerical experiments, as suggested for IBM \cite{peskin2002immersed}.

Meanwhile, equations \eqref{eq:sumf} and \eqref{eq:sumv} can be put into matrix-vector form:
\begin{align}
    \boldsymbol{f}_h(t) = & B\boldsymbol{F}_h(t), \label{eq:forcespread}\\
    \dot{\boldsymbol{X}}_h(t) = & h^2B^T\boldsymbol{u}_h(t), \label{eq:velconcentrate}
\end{align}
where 
\[
\boldsymbol{F}_h(t) = \begin{bmatrix}
\vdots\\
F_h(s,t) \ \Delta s\\
\vdots
\end{bmatrix}\in \RR^{n_s}, \quad 
\boldsymbol{X}_h(t) = \begin{bmatrix}
\vdots\\
X_h(s,t)\\
\vdots
\end{bmatrix}\in \RR^{n_s}, \quad \forall s\in \Gamma_h,
\]
are the discrete representations of the structure position and the Lagrangian force density. Using natural arrangement of the fluid variables, $B \in \RR^{n_u,n_s}$ can be constructed as a block matrix
\begin{equation}\label{eq:Bblock}
B = \begin{bmatrix}
B_1 & 0\\
0 & B_2
\end{bmatrix} \quad \left(
\begin{bmatrix}
B_1 & 0 & 0\\
0 & B_2 & 0\\
0 & 0 & B_3
\end{bmatrix} \quad \text{in 3D}\right).
\end{equation}
A column of each $B_k$ consists of evaluations of $\delta_r(\boldsymbol{x}-X(s,t))$ for a fixed $X(s,t)$ on grid points $\boldsymbol{x}$ that store one component of the fluid velocity variables. For example, the $i,j$-entry of $B_1$ is $\delta_r(\boldsymbol{x}_i-X(s_j,t))$, where $\boldsymbol{x}_i$ is the grid point that stores the $i$th fluid velocity in the $x$-direction. {\color{black}Note that $B_k$'s are not identical because they correspond to different Eulerian grid points $\boldsymbol{x}$. For example in Fig \ref{fig: semidiscrete}, $\boldsymbol{x}$ in $B_1$ are the points marked by filled triangles, while $\boldsymbol{x}$ in $B_2$ are marked by stars.} We also point out that $B$ is time dependent due to its dependence on $\boldsymbol{X}_h(t)$.

Lastly, by substituting \eqref{eq:forcespread} into \eqref{eq:DAEmotion} for $\boldsymbol{f}_h(t)$, we get the following DAE system which we shall refer to as the FOM
\begin{align}
    \rho\dot{\boldsymbol{u}}_h(t) =  & -N\boldsymbol{p}_h(t) + \mu G\boldsymbol{u}_h(t) + B\boldsymbol{F}_h(t), \label{eq:DAEsysMo}\\
    0 = & M\boldsymbol{u}_h(t), \label{eq:DAEsysMa}\\
    \dot{\boldsymbol{X}}_h(t) = & h^2B^T\boldsymbol{u}_h(t). \label{eq:DAEsysX}
\end{align}
In general, the number of structure variables is much less than the number of fluid variables, i.e., $n_s \ll n_u$. In fact, $\delta_r$ having compact support means only a small fraction of the Eulerian grid points are directly interacting with the structure. If one is only interested in the motion of the structure, i.e. $\boldsymbol{X}_h(t)$, solving the system \eqref{eq:DAEsysMo}-\eqref{eq:DAEsysX} becomes a {\color{black}reduced-order} problem \cite{freund2000krylov}, where $\boldsymbol{F}_h(t)$ is the low-dimensional input and $\dot{\boldsymbol{X}}_h(t)$ is the low-dimensional output.

\section{Reduced-order Model}\label{sec:reduce}

To construct our ROM, we start by transforming the DAE system \eqref{eq:DAEsysMo}-\eqref{eq:DAEsysX} to a coupled ODE system. Multiplying \eqref{eq:DAEsysMo} by $M$ to the left and using \eqref{eq:DAEsysMa}, we rewrite \eqref{eq:DAEsysMo} as:
\begin{equation}
    0 = -MN\boldsymbol{p}_h(t)+\mu MG\boldsymbol{u}_h(t) + MB\boldsymbol{F}_h(t).
\end{equation}
Assuming $MN$ is nonsingular, it follows that
\begin{equation}\label{eq:psol}
    \boldsymbol{p}_h(t) = (MN)^{-1}\brac{\mu MG\boldsymbol{u}_h(t) + MB\boldsymbol{F}_h(t)}.
\end{equation}
Substituting \eqref{eq:psol} into \eqref{eq:DAEsysMo} for $\boldsymbol{p}_h(t)$, one gets
\begin{equation}
    \rho\dot{\boldsymbol{u}}_h(t) = \mu QG\boldsymbol{u}_h(t) + QB\boldsymbol{F}_h(t), 
\end{equation}
where 
\begin{equation}
    Q = I-N(MN)^{-1}M  
\end{equation} 
is an oblique projection. {\color{black} It is {\color{black}worth emphasizing} here that the discrete gradient operator $N$ and the discrete divergence operator $M$ are adjoint of each other with different dimensions.}

An ODE system is then obtained from \eqref{eq:DAEsysMo},
\begin{align}
    \dot{\boldsymbol{u}}_h(t) = & \frac{\mu}{\rho} QG\boldsymbol{u}_h(t) + \frac{1}{\rho}QB\boldsymbol{F}_h(t), \label{eq:ODEsysMo}\\
    \dot{\boldsymbol{X}}_h(t) = & h^2B^T\boldsymbol{u}_h(t). \label{eq:ODEsysX}
\end{align}
We assume $\boldsymbol{u}_h(0) = \boldsymbol{0}$ in the rest of this section. Nonzero initial values can be handled by linear superposition
\begin{equation}\label{eq:linsuppos}
    \boldsymbol{u}_h(t) = \boldsymbol{u}_h^{(0)}(t) + \boldsymbol{u}_h^{(1)}(t), 
\end{equation}
in which $\boldsymbol{u}_h^{(0)}(0) = \boldsymbol{0}$ and $\boldsymbol{u}_h^{(1)}(0) = \boldsymbol{u}_h(0)$. Then one can decompose \eqref{eq:ODEsysMo}-\eqref{eq:ODEsysX} to
\begin{align}
    \dot{\boldsymbol{u}}^{(0)}_h(t) = & \frac{\mu}{\rho} QG\boldsymbol{u}^{(0)}_h(t) + \frac{1}{\rho}QB\boldsymbol{F}_h(t), \label{eq:ODEsysMo0}\\
    \dot{\boldsymbol{u}}^{(1)}_h(t) = & \frac{\mu}{\rho} QG\boldsymbol{u}^{(1)}_h(t),  \label{eq:ODEsysMo1} \\
    \dot{\boldsymbol{X}}_h(t) = & h^2B^T(\boldsymbol{u}_h^{(0)}(t) + \boldsymbol{u}_h^{(1)}(t)). \label{eq:ODEsysXdec}
\end{align}
The dynamics of $\boldsymbol{u}_h^{(1)}(t)$, which has nonzero initial value, is described by a first-order linear ODE, without interactions with the immersed structure. It can be solved separately in advance, or in some cases, it can be resolved analytically. 

We consider a general Galerkin projection of \eqref{eq:ODEsysMo}, motivated by its success in reduced-order problems \cite{freund2000krylov,Bai2002}. More specifically, we seek $\tilde{\boldsymbol{u}}_h(t)$ in a subspace, spanned by the columns of a {\color{black} tall} matrix $V$, as an approximation for $\boldsymbol{u}_h(t)$, such that for any $\boldsymbol{w}(t)$ in a test 
space, spanned by the columns of a {\color{black} tall} matrix $W$, we have
\begin{equation}
    \brac{\dot{\tilde{\boldsymbol{u}}}_h(t) - \frac{\mu}{\rho} QG\tilde{\boldsymbol{u}}_h(t) - \frac{1}{\rho}QB\boldsymbol{F}_h(t), \boldsymbol{w}(t)} = 0.
\end{equation}
{\color{black}Note that the subspaces are not necessarily fixed, which means that the matrices $V$ and $W$ are generally time-dependent. This point will be  addressed in Section 3.2.}

In a {\color{black}matrix}-vector form, the approximate solution is written as
\begin{equation}
    \tilde{\boldsymbol{u}}_h(t) = V\boldsymbol{z}(t).
\end{equation}
Then the Galerkin projection yields a reduced-order equation
\begin{equation}
    W^T\dot{(V\boldsymbol{z})}(t) = \frac{\mu}{\rho}W^TQGV\boldsymbol{z}(t) +\frac{1}{\rho}W^TQB\boldsymbol{F}_h(t).
\end{equation}

Thus we obtain an ROM of \eqref{eq:ODEsysMo} - \eqref{eq:ODEsysX}:
\begin{align}
    \dot{\boldsymbol{z}}(t) = & M_0^{-1}M_1\boldsymbol{z}(t) + M_0^{-1}M_2\boldsymbol{F}_h(t), \label{eq:RedMo}\\
    \dot{\boldsymbol{X}}_h(t) \approx & h^2B^TV\boldsymbol{z}(t), \label{eq:RedX}
\end{align}
where the matrices are given by,
\begin{equation}\label{eq:coeffMat}
    M_0 = W^TV, \quad M_1 = \frac{\mu}{\rho}W^TQGV-W^T\dot{V}, \quad M_2 = \frac{1}{\rho}W^TQB,
\end{equation}
assuming $M_0$ is nonsingular. The computation of the coefficient matrices $M_0$, $M_1$ and $M_2$ depends on $V$ and $W$. In the rest of this section we first discuss our choice for the subspaces $V$, $W$ and their properties. Then we demonstrate how an interpolation procedure can help accelerate the computation of the coefficients by exploiting the connection between the matrix entries and the Green's function.

\subsection{Subspace Selection}
We propose the following choice of $V$ and $W$,
\begin{equation}\label{eq: choicesubspaces}
    V = QB, \quad W = B.
\end{equation}
For later reference, note that both subspaces vary in time. The resulting coefficient matrices are given by
\begin{equation}\label{eq:RedMats}
    M_0 = B^TQB, \quad M_1 = \frac{\mu}{\rho}B^TQGQB-B^TQ\dot{B}, \quad M_2 = \frac{1}{\rho}B^TQB = \frac{1}{\rho}M_0.
\end{equation}

{\color{black}In principle, one can use higher dimensional Krylov subspaces ($V$ and $W$ with more columns), followed by Lanczos orthogonalization algorithms \cite{Bai2002,freund2000krylov,ma2019coarse}, to improve the accuracy of the ROM. Specifically, we shall see in the following discussion that our choice satisfies two interpolation conditions. Higher dimensional Krylov subspaces are able to interpolate the transfer function more accurately by enforcing more interpolation conditions, {\color{black}but at the expense of a reduced computational speedup.} From the numerical tests, our observation is that the current subspaces achieve a good balance between accuracy and efficiency.}

\subsubsection{Transfer function approximation}
{\color{black}We first show the accuracy property of our choice of subspaces. This can be understood by solving the linear ODE \eqref{eq:ODEsysMo} analytically for $\boldsymbol{u}_h$, which yields,
\begin{equation}\label{eq:u_ana}
    \boldsymbol{u}_h(t) = \frac{1}{\rho}\int_0^t \exp\big[\frac{\mu}{\rho}(t-\tau)QG\big]QB\boldsymbol{F}_h(\tau)\, d\tau.
\end{equation}
Note that we assume zero initial condition as discussed before. Plugging \eqref{eq:u_ana} into \eqref{eq:ODEsysX} gives
\begin{equation}
    \dot{\boldsymbol{X}}_h(t) = \int_0^t \phi(t-\tau)\boldsymbol{F}_h(\tau)\, d\tau,
\end{equation}}
where $\phi(t)$ denotes the transfer function,
\begin{equation}\label{eq:fulltf}
    \phi(t) = \frac{h^2}{\rho}B^T\exp\big[\frac{\mu}{\rho}tQG\big]QB.
\end{equation}
A similar calculation for the reduced-order model \eqref{eq:RedMo} - \eqref{eq:RedX} shows that: 
\begin{equation}\label{eq:approxXdot}
    \dot{\boldsymbol{X}}_h(t) \approx \int_0^t \phi_{red}(t-\tau)\boldsymbol{F}_h(\tau)\, d\tau,
\end{equation}
where the transfer function $\phi_{red}(t)$ of the reduced-order model is given by,
\begin{align}
    \phi_{red}(t) = & h^2B^TV\exp[tM_0^{-1}M_1]M_0^{-1}M_2, \nonumber \\
           = & \frac{h^2}{\rho}B^TQB\exp\left[t(B^TQB)^{-1}\brac{\frac{\mu}{\rho}B^TQGQB-B^TQ\dot{B}}\right].\label{eq:Redtf}
\end{align}
$\phi_{red}(t)$ is expected to approximate $\phi(t)$ in the sense that,
\begin{align}
    \phi_{red}(0) = \phi(0),\label{eq:tfzero}\\
    \dot{\phi}_{red}(0) = \dot{\phi}(0). \label{eq:tfdiffzero}
\end{align}
The equality \eqref{eq:tfzero} follows immediately from evaluating \eqref{eq:fulltf} and \eqref{eq:Redtf} at $t=0$. Differentiating \eqref{eq:fulltf} and \eqref{eq:Redtf} at $t=0$ yields
\begin{equation}
    \dot{\phi}_{red}(0) = \frac{\mu h^2}{\rho^2}B^TQB(B^TQB)^{-1}B^TQGQB = \frac{\mu h^2}{\rho^2}B^TQGQB = \dot{\phi}(0).
\end{equation}
In the above calculation, we have treated $B$ as a constant matrix. The reason is that we are only concerned with a small time interval $[0,t]$, typically with the size of one time step. In numerical simulations, the matrix $B$ is usually treated as constant when advancing one time step.
\subsubsection{Enforcing incompressibility}
Another essential property of the full model is the incompressibility of the fluid. Recall that $M$ is the discrete divergence operator. The approximate fluid solution,
\begin{equation}\label{eq:fluidappr}
\tilde{\boldsymbol{u}}_h(t) = V\boldsymbol{z}(t),
\end{equation}
is incompressible if
\begin{equation}\label{eq:RedDivfree}
    M\tilde{\boldsymbol{u}}_h(t) = MV\boldsymbol{z}(t) = 0.
\end{equation}
A quick calculation verifies that our choice of $V=QB$ satisfies this constraint:
\begin{equation}
    MV = MQB = M(I-N(MN)^{-1}M)B = (M-MN(MN)^{-1}M){\color{black}B = 0.}
\end{equation}
Therefore, the incompressibility property is preserved in the ROM. 
\subsubsection{Lyapunov Stability}
The ROM also preserves Lyapunov stability of the FOM with our choice of subspaces. We first show the stability of the FOM. We assume the discrete Lagrangian force density $\boldsymbol{F}_h$ is given by an energy functional $W(\boldsymbol{X}_h)$ of the structure configuration, i.e.,
\begin{equation}
\boldsymbol{F}_h(t) = -\nabla_{X_h} W(\boldsymbol{X}_h(t)).
\end{equation}
We also assume the discrete gradient and divergence operators satisfy
\begin{equation}\label{eq:graddivsym}
M=N^T,
\end{equation}
such that
\begin{equation}\label{eq:Qorthproj}
Q = I-\Sigma = I-N(N^TN)^{-1}N^T
\end{equation}
is an orthogonal projection. We now define the following Lyapunov functional for the FOM consisting of the kinetic and the elastic energy,
\begin{equation}
{\color{black}V(\boldsymbol{u}_h(t),\boldsymbol{X}_h(t)) = \frac{1}{2}\boldsymbol{u}_h(t)^TQ \boldsymbol{u}_h(t) + \frac{1}{h^3\rho}W(\boldsymbol{X}_h(t)).}
\end{equation}
We have $V(\boldsymbol{u}_h(t),\boldsymbol{X}_h(t))\geq 0$ because $Q$, as a projection, is positive semidefinite with eigenvalues $0$ or $1$. In particular, notice that $Q^2=Q$ and $Q^T=Q$. {\color{black}In addition, the divergence-free condition implies that $Q\boldsymbol{u}_h(t)=\boldsymbol{u}_h(t)$. A direct calculation shows that
\begin{align}
\dot{V}(\boldsymbol{u}_h(t),\boldsymbol{X}_h(t)) = & \boldsymbol{u}_h(t)^TQ^2(\frac{\mu}{\rho}G\boldsymbol{u}_h(t) + \frac{1}{\rho}B\boldsymbol{F}_h(t)) -\frac{1}{\rho}\boldsymbol{u}_h(t)^TB\boldsymbol{F}_h(t)) \nonumber \\
= & \frac{\mu}{\rho}\boldsymbol{u}_h(t)^TQG\boldsymbol{u}_h(t)+\frac{1}{\rho}\boldsymbol{u}_h(t)^TQB\boldsymbol{F}_h(t))-\frac{1}{\rho}\boldsymbol{u}_h(t)^TB\boldsymbol{F}_h(t)) \nonumber\\
= & \frac{\mu}{\rho}\boldsymbol{u}_h(t)^TG\boldsymbol{u}_h(t)\leq 0,
\end{align}}
since the discrete Laplace operator $G$ is negative semidefinite. This implies the Lyapunov stability of the FOM.

The Lyapunov functional for the ROM is defined as follows
\begin{equation}\label{redLyap}
V_r(\boldsymbol{X}_h(t)) = \frac{h^2}{\rho}W(\boldsymbol{X}_h(t))+\frac{1}{2}\dot{\boldsymbol{X}}_h(t)^T(B^TQB)^{-1}\dot{\boldsymbol{X}}_h(t).
\end{equation}
It is now clear that $V_r(\boldsymbol{X}_h(t))\geq 0$ holds for all $t$ since $B^TQB$ is positive semidefinite.

To prove $\dot{V}_r(\boldsymbol{X}_h(t)) \leq 0$, we start by rewriting the ROM \eqref{eq:RedMo} - \eqref{eq:RedX} as a second-order ODE of $\boldsymbol{X}_h$. {\color{black}Note that \eqref{eq:RedX} and \eqref{eq: choicesubspaces} imply  $\boldsymbol{z}(t)=h^{-2}(B^TQB)^{-1}\dot{\boldsymbol{X}}_h(t)$ in the ROM.} Using the symmetry of $Q$, one has,
{\color{black}\begin{align}
\ddot{\boldsymbol{X}}_h = & h^2[\dot{(B^TQB)}\boldsymbol{z}+B^TQB\dot{\boldsymbol{z}}(t)] \nonumber\\
= & h^2(B^TQ\dot{B}+\frac{\mu}{\rho}B^TQGQB)\boldsymbol{z}+\frac{h^2}{\rho}B^TQB\boldsymbol{F}_h\brac{\boldsymbol{X}_h(t)} \nonumber \\
= & (B^TQ\dot{B}+\frac{\mu}{\rho}B^TQGQB)(B^TQB)^{-1}\dot{\boldsymbol{X}}_h+\frac{h^2}{\rho}B^TQB\boldsymbol{F}_h\brac{\boldsymbol{X}_h(t)}. \label{eq:Red2ndODE}
\end{align}}
Then the following calculation shows that $V_r$ is nonincreasing,
\begin{align}\label{redLyapdot}
\dot{V}_r(\boldsymbol{X}_h(t)) = & -\frac{h^2}{\rho}\dot{\boldsymbol{X}}_h^T\boldsymbol{F}_h +\dot{\boldsymbol{X}}_h(B^TQB)^{-1}\ddot{\boldsymbol{X}}_h +\frac{1}{2}\dot{\boldsymbol{X}}_h\dot{(B^TQB)^{-1}}\dot{\boldsymbol{X}}_h \nonumber \\
= &  \frac{\mu}{\rho}\dot{\boldsymbol{X}}_h^T(B^TQB)^{-1}(B^TQGQB)(B^TQB)^{-1}\dot{\boldsymbol{X}}_h 
 -\frac{h^2}{\rho}\dot{\boldsymbol{X}}_h^T\boldsymbol{F}_h + \frac{h^2}{\rho}\dot{\boldsymbol{X}}_h^T\boldsymbol{F}_h \nonumber \\
&+\dot{\boldsymbol{X}}_h^T(B^TQB)^{-1}(B^TQ\dot{B})(B^TQB)^{-1}\dot{\boldsymbol{X}}_h \nonumber \\  
& -\dot{\boldsymbol{X}}_h^T(B^TQB)^{-1}(B^TQ\dot{B})(B^TQB)^{-1}\dot{\boldsymbol{X}}_h \nonumber \\
= &  \frac{\mu}{\rho}\dot{\boldsymbol{X}}_h^T(B^TQB)^{-1}(B^TQGQB)(B^TQB)^{-1}\dot{\boldsymbol{X}}_h \nonumber \\
= & \frac{\mu}{\rho}\boldsymbol{Y}^TG\boldsymbol{Y} \leq 0,
\end{align}
where we have defined $\boldsymbol{Y}:=QB(B^TQB)^{-1}\dot{\boldsymbol{X}}_h$. The last inequality holds because $G$ is negative semidefinite.

\subsection{Computing the coefficients using interpolation}
Because the coefficient matrices in the ROM are in principle time dependent, they should be updated frequently during simulation. {\color{black}Direct matrix multiplication for this purpose is time consuming since $Q$ is a dense matrix in $\RR^{2n_u,2n_u}$. For example, computing $M_1$ in \eqref{eq:RedMats} has complexity $O(n_u^2)$.} In the rest of this section, we propose a computationally cheaper approach using interpolation to approximate the coefficients.

We first approximate $\dot{B}$ by
\begin{equation}
    \dot{B}(t) \approx \frac{1}{\Delta t}\brac{B(t)-B(t-\Delta t)}.
\end{equation}
{\color{black}One could consider a higher order discretizations for $\dot{B}$ so that a method of order higher than one in time can be used to solve the ROM. Ultimately, this is a trade-off between accuracy and the offline interpolation efficiency. From the numerical tests, we will see such first order approximation of $\dot{B}$ together with the forward Euler's method provides acceptable accuracy for the numerical tests compared to the FOM.}  

The other observation is that the matrix $M_1$ is then approximated by
\begin{equation}\label{eq:Bdotappr}
M_1 \approx \frac{\mu}{\rho}B(t)^TQGQB(t)-\frac{1}{\Delta t}B(t)^TQB(t) + \frac{1}{\Delta t}B(t)^TQB(t-\Delta t).
\end{equation}
Together with $M_0 = B(t)^TQB(t)$, the following three matrices are needed for building our ROM
\begin{equation}\label{eq:matneeded}
B(t)^TQB(t), \quad B(t)^TQGQB(t), \quad B(t)^TQB(t-\Delta t).
\end{equation}
Since $Q$ and $G$ are constant matrices, the $i,j$-entry of any of the above matrices at time $t$ is determined by the $i$th row of $B(t)^T$ and the $j$th column of $B(t)$ (or $B(t-\Delta  t)$). Recall that each column of $B(t)$ (or row of $B(t)^T$) represents a smoothed delta function associated with a structure point. Suppose the $i$th column of $B(t)^T$ is associated with the Lagrangian grid point $X_l\in\RR^d$ ($d=2,3$) and the $j$th column of $B(t)^T$ (or $B(t-dt)$) is associated with $X_r\in\RR^d$. Given the prescribed function $\delta_r$ and a fixed Eulerian grid, the $i,j$-entry of a coefficient matrix is uniquely determined by $X_l$ and $X_r$, which can be viewed as a function from $\RR^{2d}$ to $\RR$. It is then natural to sample such functions before the simulation starts. As the simulation runs, coefficient matrices are updated by interpolation using precomputed samples. {\color{black}In this work, linear interpolation is used. Because each entry of the $2n_s$-by-$2n_s$ coefficient matrix is obtained by evaluating a precomputed linear function, the complexity is typically $O(n_s^2)$, which is much smaller than the complexity of direct matrix multiplications $O(n_u^2)$, given $n_s \ll n_u$.}

\begin{figure}[htp]
    \centering    
    \includegraphics[width = .9\textwidth]{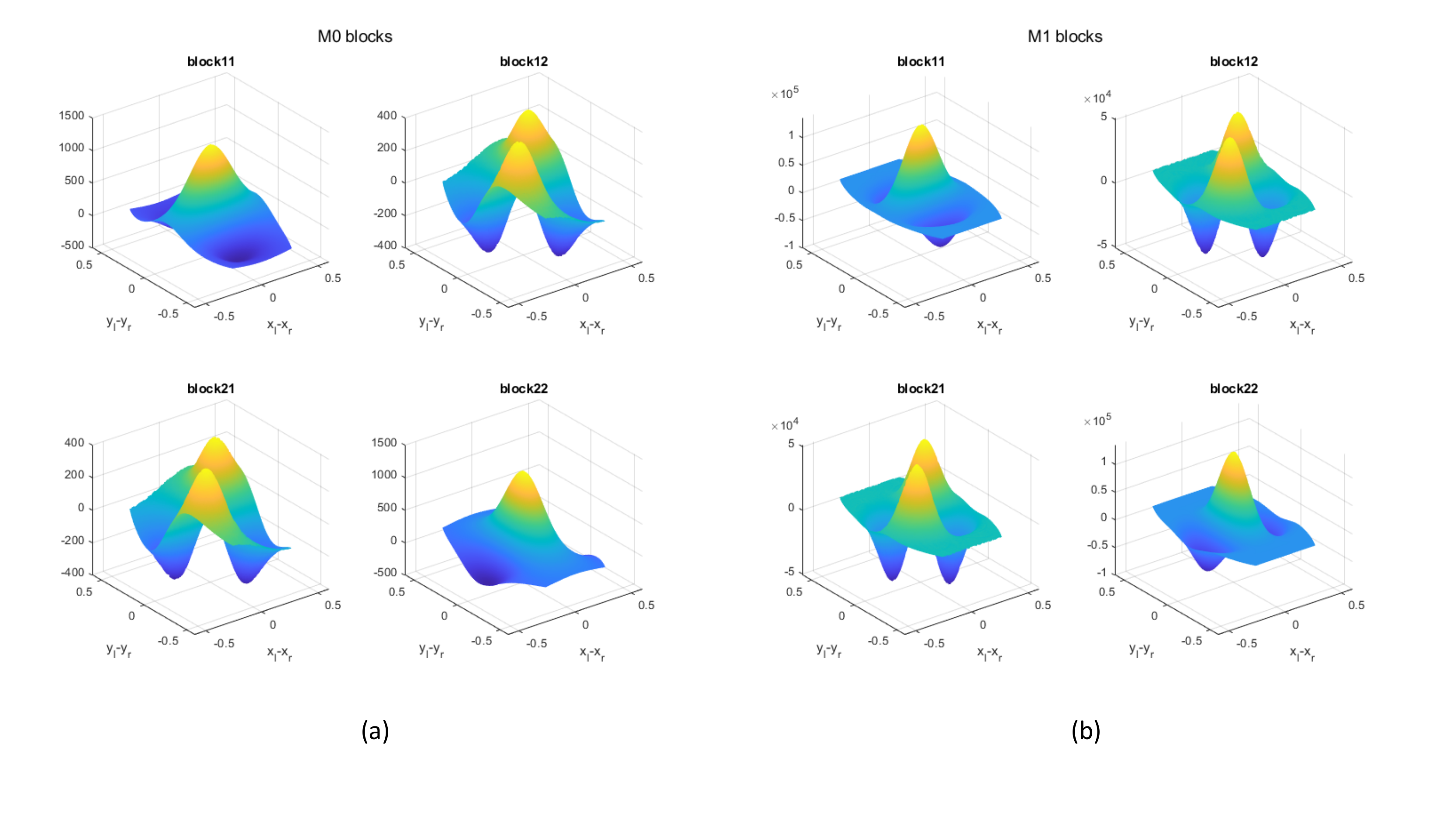}
    \caption{Surface plots of entries of (a)$M_0$ and (b)$M_1$ against $X_l-X_r = (x_l-x_r,y_l-y_r)$ in 2D case. Both matrices are 2-by-2 block matrices where each block corresponds to a function from $\RR^4$ to $\RR$. Large numbers of $(X_l,X_r) = (x_l, y_l, x_r,y_r)$ pairs are sampled so that many of them correspond to the same difference $X_l-X_r$. Then corresponding entries of $M_0$ and $M_1$ blocks are plotted against $X_l-X_r$. In each plot, we observe a single surface, indicating no multiple values. Therefore these $\RR^4$ functions can be considered as functions in $\RR^2$ of $X_l-X_r$. }
    \label{fig: fcn2d}
\end{figure}

Next, motivated by our numerical experiments illustrated in Fig \ref{fig: fcn2d}, we show that the interpolated $2d$-dimensional functions of $X_l$ and $X_r$ can be well approximated by $d$-dimensional functions of $X_l-X_r$, 
{\color{black} i.e., the relative position of the two points.} Such low-dimensional approximation significantly reduces the number of samples needed for more accurate interpolations. Hence the sampling process can also be accelerated. 
\begin{figure}[htp]
    \centering    
    \includegraphics[width = .9\textwidth]{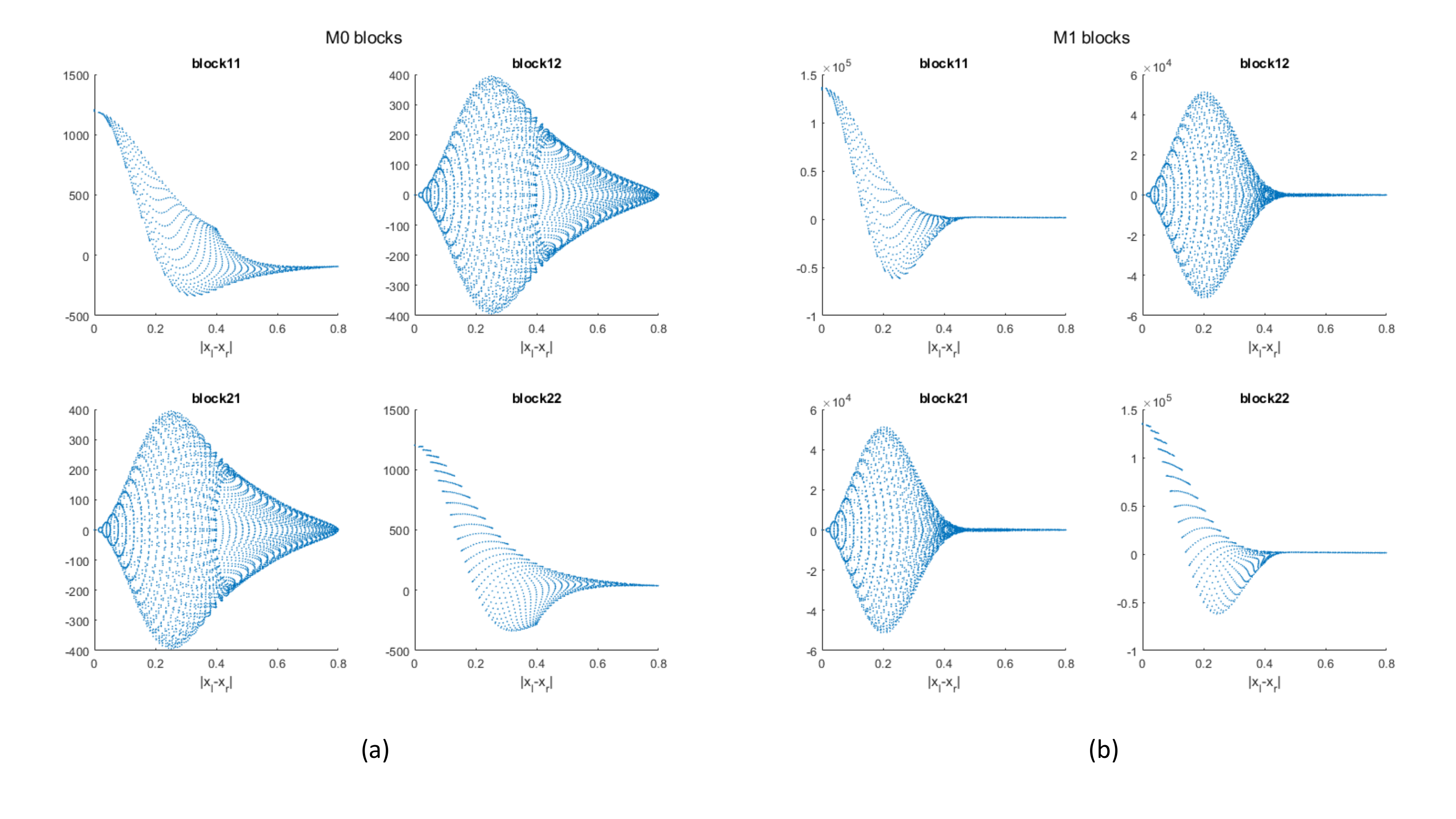}
    \caption{Scatter plots of entries of (a)$M_0$ and (b)$M_1$ against $|X_l-X_r|$ in 2D case. Multiple values exist, indicating these $\RR^4$ functions may not be considered as functions in $\RR$ of $|X_l-X_r|$. }
    \label{fig: fcn1d}
\end{figure}

Here we provide justifications of this approach by making connections to the Green's functions. 
Recall that $Q = I-N(MN)^{-1}M$, where $N$ and $M$ are discrete gradient and divergence operators. Therefore, each entry of the matrix $B(t)^TQB(t)$ or $B(t)^TQB(t-\Delta t)$ is a numerical approximation of the integral
\begin{align}
    I_0 = & \int_{\Omega} \delta_r(\boldsymbol{x}-X_l)(\delta_{ij}-\partial_i \Delta^{-1}\partial_j) \delta_r(\boldsymbol{x}-X_r)\,d\boldsymbol{x} \nonumber \\
    = & \underbrace{\delta_{ij}\int_{\Omega} \delta_r(\boldsymbol{x}-X_l)\delta_r(\boldsymbol{x}-X_r)\,d\boldsymbol{x}}_{I_1} - \underbrace{\int_{\Omega} \delta_r(\boldsymbol{x}-X_l)\partial_i \Delta^{-1}\partial_j \delta_r(\boldsymbol{x}-X_r)\,d\boldsymbol{x}}_{I_2}, \label{eq:M0int}
\end{align}
where $i,j=1,\cdots,d$ and $d=2$ or $3$. $\delta_{ij}$ is the Kronecker delta function. $I_1$ only depends on $|X_l-X_r|$ due to our choice of $\delta_r$. For $I_2$, we assume $X_l$ and $X_r$ are far from the boundary of $\Omega$ so the Green's function $G(\boldsymbol{x},\boldsymbol{y})$ can be applied. Considering the limiting case of $r\to 0$, i.e.,  $\delta_r \to \delta$, as $r \to 0,$ and we arrive at,
\begin{equation}
    \lim_{r\to 0}I_2 = -\partial_{x_i} \partial_{y_j} G(\boldsymbol{x},\boldsymbol{y})\big\vert_{\boldsymbol{x}=X_l,\boldsymbol{y}=X_r},
\end{equation}
which depends only on $X_l-X_r$.

Similarly, each entry of the matrix $B(t)^TQGQB(t)$ is a numerical approximation of the following integral
\begin{align}
    J_0 = & \int_{\Omega} \delta_r(\boldsymbol{x}-X_l)(\delta_{ij}-\partial_i \Delta^{-1}\partial_j)\Delta(1-\partial_j \Delta^{-1}\partial_j) \delta_r(\boldsymbol{x}-X_r)\,d\boldsymbol{x} \nonumber \\
    = & \underbrace{\int_{\Omega} \delta_r(\boldsymbol{x}-X_l)(\delta_{ij}-\partial_i \Delta^{-1}\partial_j)\Delta \delta_r(\boldsymbol{x}-X_r)\,d\boldsymbol{x}}_{J_1} 
      - \underbrace{\int_{\Omega} \delta_r(\boldsymbol{x}-X_l)(\delta_{ij}-\partial_i \Delta^{-1}\partial_j)\Delta\partial_j \Delta^{-1}\partial_j \delta_r(\boldsymbol{x}-X_r)\,d\boldsymbol{x}.}_{J_2}  \nonumber \\
    = & \underbrace{\delta_{ij}\int_{\Omega} \delta_r(\boldsymbol{x}-X_l)\Delta \delta_r(\boldsymbol{x}-X_r)\,d\boldsymbol{x}}_{J_{11}} 
        - \underbrace{\int_{\Omega} \delta_r(\boldsymbol{x}-X_l)\partial_i \Delta^{-1}\partial_j\Delta \delta_r(\boldsymbol{x}-X_r)\,d\boldsymbol{x}}_{J_{12}} \nonumber \\
        & - \underbrace{\delta_{ij}\int_{\Omega} \delta_r(\boldsymbol{x}-X_l)\partial_j^2 \delta_r(\boldsymbol{x}-X_r)\,d\boldsymbol{x}}_{J_{21}} 
        + \underbrace{\int_{\Omega} \delta_r(\boldsymbol{x}-X_l)\partial_i \Delta^{-1}\partial_j\Delta\partial_j \Delta^{-1}\partial_j \delta_r(\boldsymbol{x}-X_r)\,d\boldsymbol{x}}_{J_{22}},     \label{eq:M1int}
\end{align}
where $J_{11}$ and $J_{21}$ depend on $X_l-X_r$ due to our choice of $\delta_r$. For $J_{12}$ and $J_{22}$ we make the same assumptions as for $I_2$ and consider the limiting case. We obtain similar results
\begin{align}
    \lim_{r\to 0}J_{12} = & -\partial_{x_i} \partial_{y_j}^3 G(\boldsymbol{x},\boldsymbol{y})\big\vert_{\boldsymbol{x}=X_l,\boldsymbol{y}=X_r},\\
    \lim_{r\to 0}J_{22} = & -\partial_{x_i} \Delta_y\partial_{y_j} G(\boldsymbol{x},\boldsymbol{y})\big\vert_{\boldsymbol{x}=X_l,\boldsymbol{y}=X_r}.
\end{align} So both terms  depend only on $X_l-X_r$. However, the integrals $I_0$ and $J_0$ may not be further reduced to functions of $|X_l-X_r|$, as suggested by our numerical experiments, see Fig \ref{fig: fcn1d}.

\section{Numerical results}\label{sec:numres}
In this section, we present three numerical examples to demonstrate the accuracy and speedup offered by our ROM. {\color{black}The finite difference method is used for both FOM and ROM. For the temporal discretization, we use the forward Euler method.}

\subsection{Oscillation of an elliptical membrane}
\begin{figure}[ht]
    \centering
    \includegraphics[width=1\textwidth]{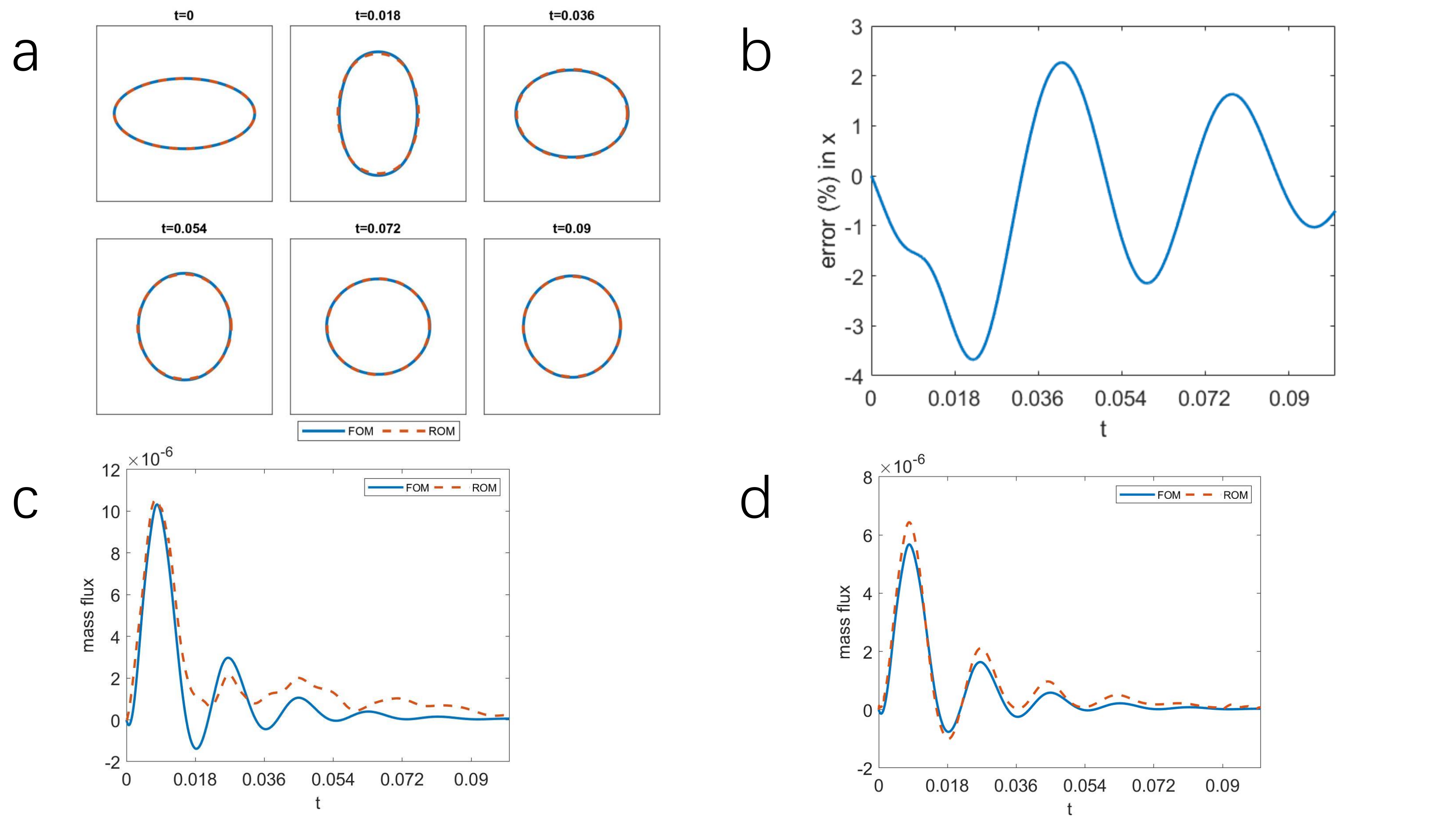}
    \caption{(a) Comparison between FOM and ROM of the elliptical membrane profiles at different times. (b) Relative difference between the x-coordinates of a reference structure point simulated by the FOM and ROM. (c, d) Evolution of mass flux across the membrane with fluid grid size $h=1/8$ (c) and $h=1/16$ (d).}
    \label{fig: memsnaps}
\end{figure}

We consider the oscillations of a pressurized fiber. {\color{black}Initially, the stretched elastic fiber resides in the center of a resting fluid. The semi-major and semi-minor axes of the fiber are 0.4 and 0.2 $\mu m$, respectively. The fluid domain is $4 \mu m\times 4\mu m$ with periodic boundary conditions on all edges. Fluid density and viscosity are chosen so that the Reynolds number is $0.01$.} The body force in this example is generated by an elastic energy functional \cite{peskin2002immersed},
\begin{equation}\label{eq:ex1energy}
    E = \int_\Gamma \varepsilon\left(\left\vert\pd{\boldsymbol{X}}{s}\right\vert\right)\,ds,
\end{equation}
where $\varepsilon$ is the local energy given by
\begin{equation}\label{eq:ex1energyfunc}
    \varepsilon(x) = \frac{\sigma}{2}(x-L)^2,
\end{equation}
which corresponds to an elastic fiber having a "spring constant" $\sigma$ and an equilibrium state where the elastic strain $\vert\partial\boldsymbol{X}/\partial s\vert \equiv L$. The force in \eqref{eq:forceDensity} is then expressed as  
\begin{equation}\label{eq:ex1force}
    F = -\nabla_XE = \sigma\pd{}{s}\left( \pd{\boldsymbol{X}}{s}\left(1-\frac{L}{\left\vert \pd{\boldsymbol{X}}{s} \right\vert} \right)  \right).
\end{equation}
Since the fluid in the interior of the membrane is confined, the membrane will oscillate and eventually settle into a circular state. Membrane configurations simulated by the FOM and ROM are compared at different times (Fig \ref{fig: memsnaps}a \& b). The ROM simulation captures almost the same equilibrium state as the FOM. In addition, the membrane configurations are approximated accurately during the oscillation. {\color{black} We demonstrate that the ROM preserves the incompressibility by comparing the evolution of mass flux with that of the FOM. The mass flux is calculated by integrating the velocity over the membrane surface using the trapezoidal rule. The mass flux of the ROM is in close agreement with the FOM. Both are very close to zero up to a numerical error which keeps decreasing as the grid becomes finer, as shown in Fig. \ref{fig: memsnaps}c \& d.} 

The one-step computation time of our ROM simulations with various grid sizes is compared to the one-step FOM simulation time in Table \ref{tab: memspeedup}. There is a clear increase in the speedup factor as the grid spacing decreases. With 2D flow, the {\color{black}time complexities} are $O(h^{-2})$ and $O(h^{-1})$ for the FOM and ROM simulations, respectively. {\color{black}The effect of the additional sampling cost at the beginning of the simulation is reported in Table \ref{tab: memcost}. This overhead is less than $20$ time steps of the FOM simulation. In this example, the total number of time steps is $1000$. Therefore the computational cost associated with the sampling process is negligible compared to the speedup during the simulation.}

\begin{figure}[ht]
    \centering
    \includegraphics[width=.6\textwidth]{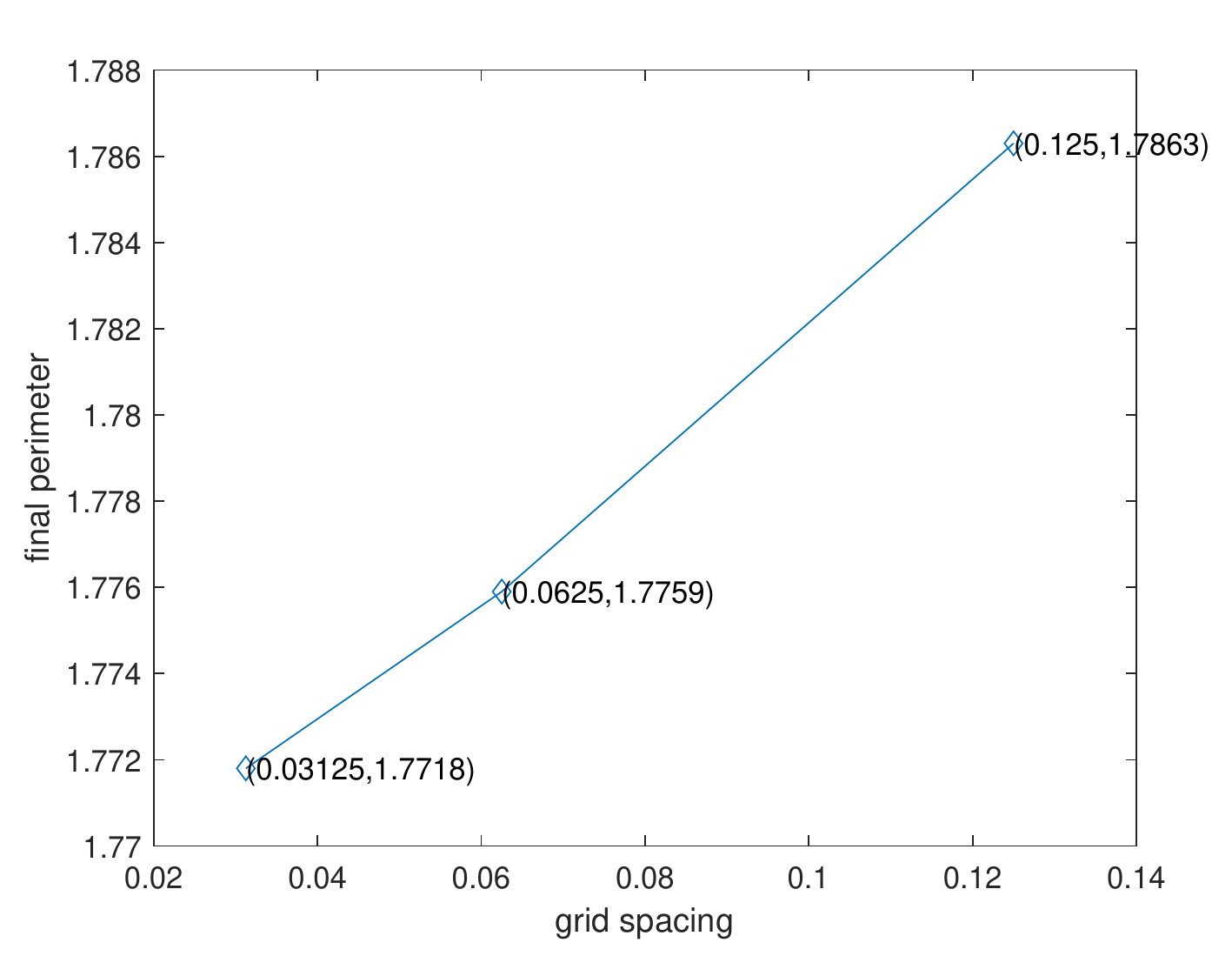}
    \caption{Convergence plot of ROM measured by final perimeter with varying grid spacings.}
    \label{fig: conv}
\end{figure}

{\color{black}We show the perimeter of the immersed structure at final time for various choices of the grid size (Fig \ref{fig: conv}). As the grid spacing reduces, the perimeter approaches an asymptotic zero-grid spacing value. We determine the order of convergence of the ROM based on these results,
\[
\ln\left(\frac{1.7759-1.7863}{1.7718-1.7759}\right) / \ln(2) = 1.3429.
\]}

\subsection{Rotation of an elliptical particle in shear flow}
\begin{figure}
    \centering
    \includegraphics[width=1\textwidth]{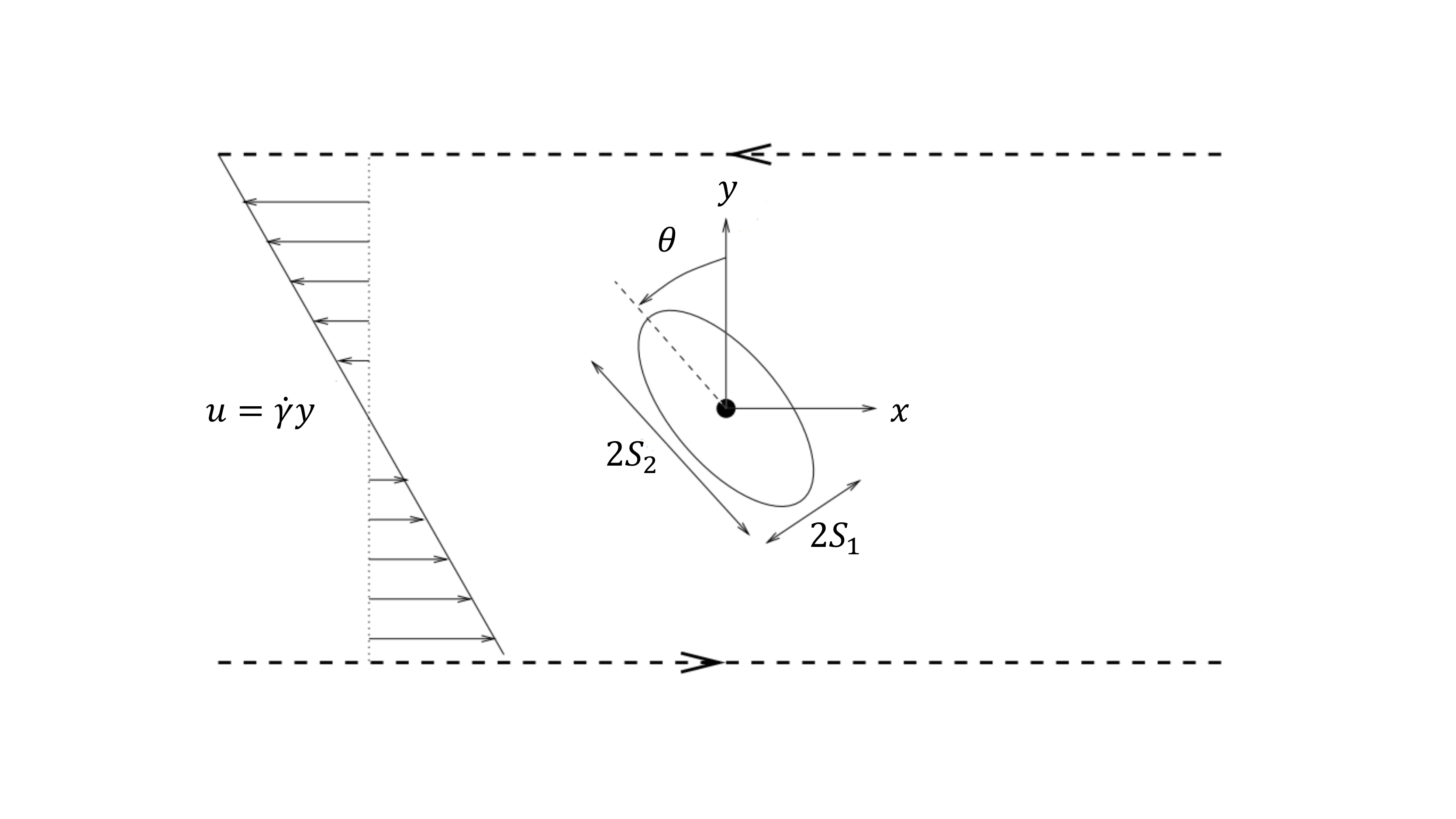}
    \caption{A rigid ellipse immersed in a shear flow.}
    \label{fig: orbdem}
\end{figure}

We study the problem of the motion of a rigid elliptical particle freely suspended in a shear flow. {\color{black}The fluid domain is $8\mu m \times 8\mu m$. The semi-minor and semi-major axes of the ellipse are $S_1=0.2 \mu m$ and $S_2=0.3 \mu m$, respectively. Initially, the ellipse is immersed in the center of a shear flow with its semi-major axis positioned along the $y$-axis. The maximum fluid velocity of the shear flow, fluid density, viscosity are chosen so that the Reynolds number is $0.01$. (Fig \ref{fig: orbdem})} 

It has been shown that the instantaneous inclination angle $\theta$ of the ellipse major axis with respect to the $y$-axis is
\begin{equation}\label{eq:orbittrue}
    \tan(\theta) = \frac{S_2}{S_1}\tan\brac{\frac{S_1S_2}{S_1^2+S_2^2}\dot{\gamma}t},
\end{equation}
where $t$ is the time variable \cite{jeffery1922motion}. 

To preserve the elliptic shape of the rigid structure, the body force in this example is generated by a discrete bending energy {\color{black}\cite{pivkin2008accurate}}. Let $\theta_i^0$ be the initial angle between the adjacent edges with the $i$-th Lagrangian grid point and $\theta_i$ be the current angle. The bending energy is given by
\begin{equation}
E_b = \sigma_b \sum_{i=1}^{n_s}(1-\cos(\theta_i-\theta_i^0)),
\end{equation}
where $n_s$ is the number of Lagrangian grid points and $\sigma_b$ is the bending coefficient. In this example, we choose $\sigma_b=2000$ to increase the stiffness. The bending force generated on each structure point is given by, 
\begin{equation}
\mathbf{F}_i = ({F_i} _x,{F_i} _y) = (-\frac{\partial E_b}{\partial x_i},-\frac{\partial E_b}{\partial y_i} ).
\end{equation}
\begin{figure}
    \centering
    \includegraphics[width=1\textwidth]{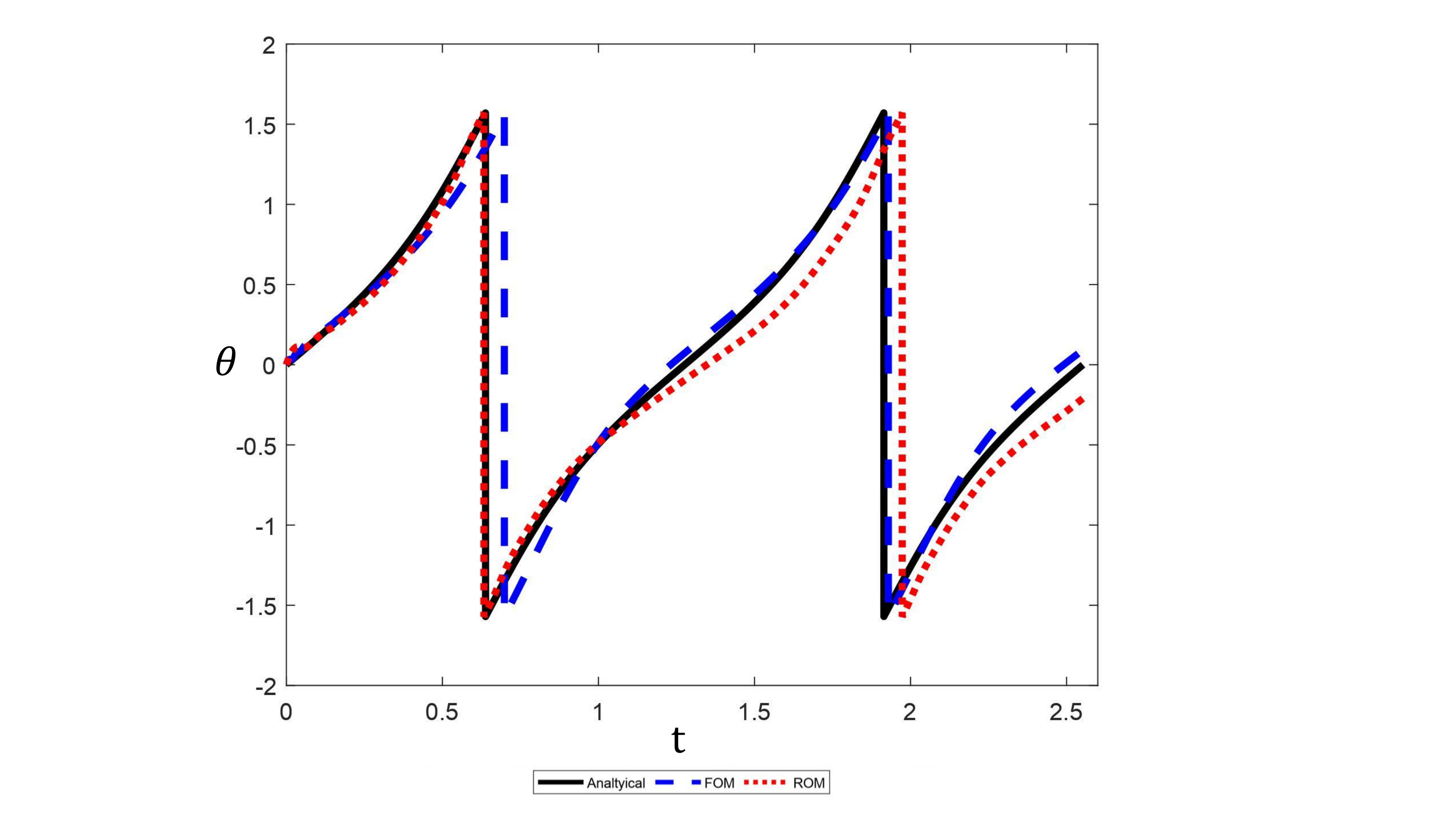}
    \caption{Ellipse rotation angles simulated by the full model and the reduced-order model compared with Jeffery’s orbit. The variation in the angle $\theta$ relative to the ellipse major axis is plotted as a function of time $t$.}
    \label{fig: orbangles}
\end{figure}

Fig \ref{fig: orbangles} shows the simulated ellipse rotation rate and the analytical result \eqref{eq:orbittrue}. The rotation rate obtained by our ROM simulation is in close agreement with both the FOM simulation and the analytical solution.Table \ref{tab: orbspeedup} shows the increase in the speedup factor as the grids become finer. Higher speedup factors are achieved for finer space grid.
\subsection{Motion of two particles in laminar flow}
\begin{figure}
    \centering
    \includegraphics[width=.8\textwidth]{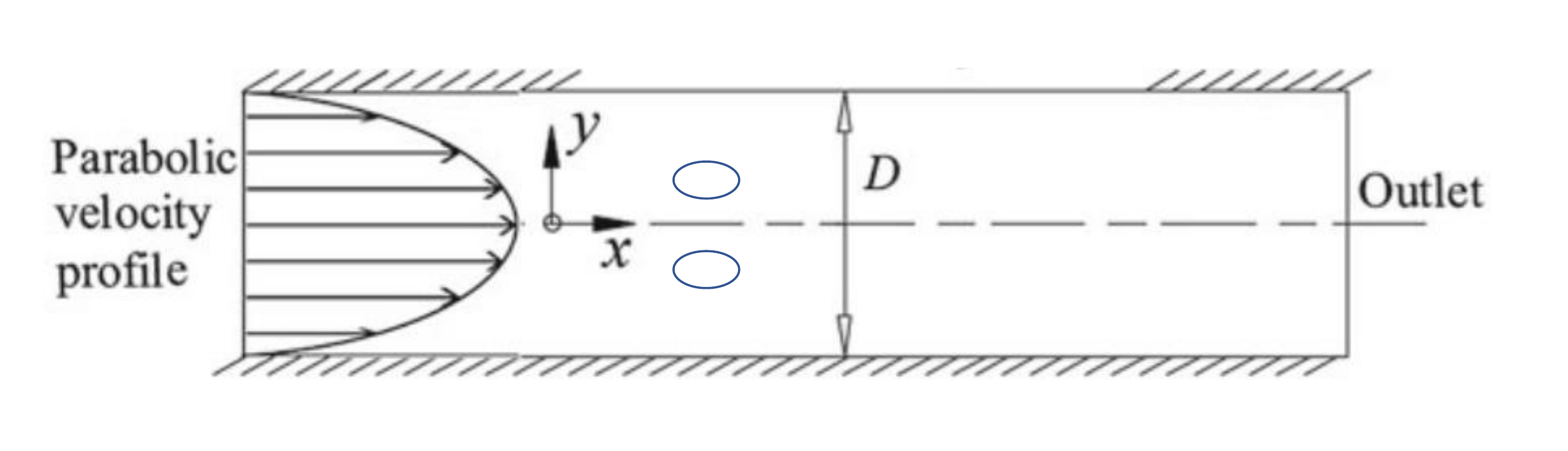}
    \caption{Two membranes interacting in a laminar channel flow.}
    \label{fig: ex3dem}
\end{figure}

{\color{black} In the last numerical test, we simulate the motion of two membranes in a $6\mu m \times 15\mu m$ channel. The fluid is initially at rest, with inlet velocity profile given by,  as depicted in Fig \ref{fig: ex3dem},
\begin{equation}
    U = U_0\left[ 1-\left(\frac{y}{D}\right)^2 \right], \quad -D\leq y\leq D.
\end{equation}
At the beginning, the two membranes of the same elliptic shape are placed with horizontal semi-major axes and the same distance $0.6 \mu m$ from its center to the x-axis. The initial semi-major axis and semi-minor axis are $0.3 \mu m$ and $0.2 \mu m$, respectively. Fluid density, viscosity, $U_0$ are chosen so that the Reynolds number is $0.01$. Nonslip conditions are applied to the top and bottom boundaries.} 

\begin{figure}[htp]
    \centering    
    \includegraphics[width = 1\textwidth]{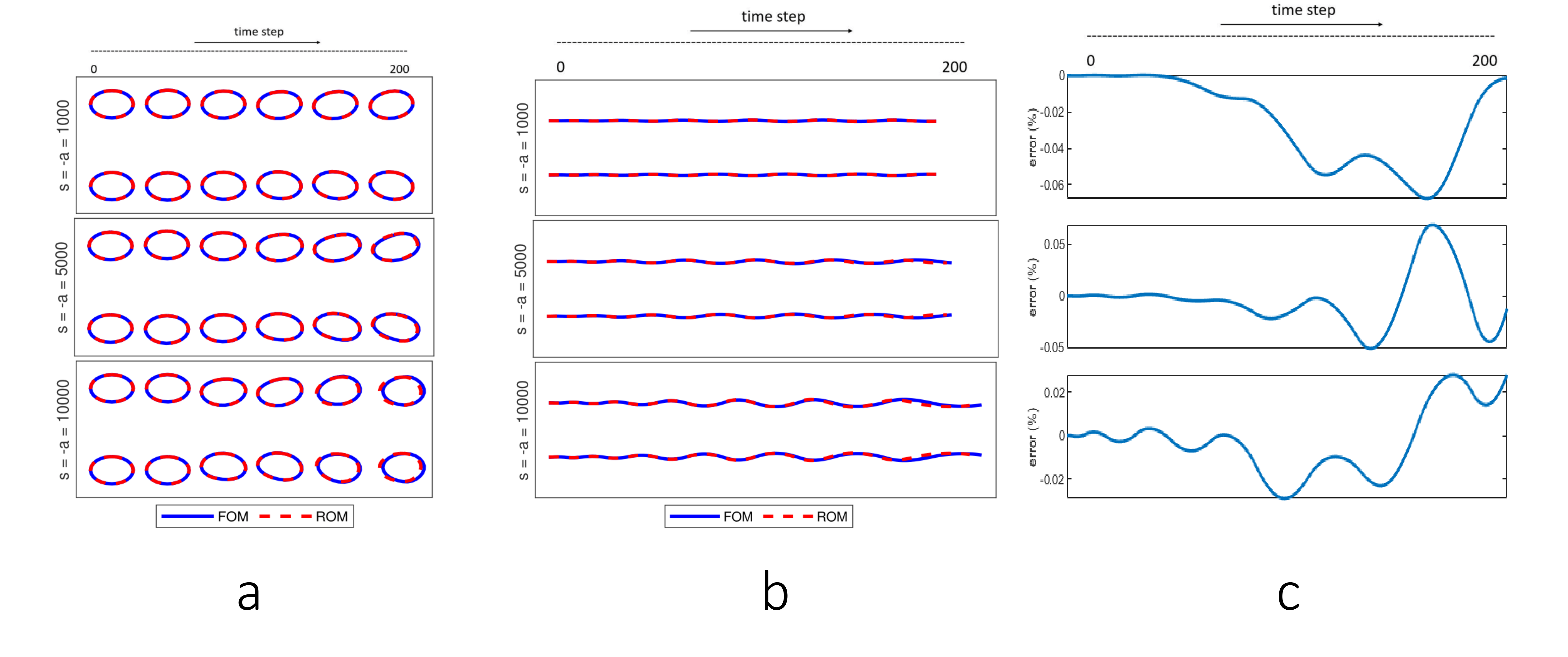}
    \caption{The increase of the attraction and repulsion forces (from top to bottom) as the magnitude of $a$ and $s$ increases. The other two parameter, $b$ and $\lambda$, are fixed for ease of comparison. (a) Snapshots of two membranes at different times. (b) Trajectories of membrane centers. (c) Relative error in x-coordinates of a reference structure point of the upper membrane, between the FOM and ROM.}
    \label{fig: tctrajsnaps}
\end{figure}

The same bending force as in the previous example is applied to both membranes to prevent significant deformation. In addition, the two membranes interact with each other through a binding force and a repulsive force given respectively by,
\begin{align}
    F_{binding} = & s(d-\lambda), \label{eq:twocell_bind} \\
    F_{repulsion} = & ad+bd^3, \label{eq:twocell_rep}
\end{align}
where $d$ is the distance between two Lagrangian nodes on different cells and $a$, $b$, $s$, $\lambda$ are parameters. These forces are developed to model the biochemical interactions between flowing melanoma tumor cells and substrate adherent polymorphonuclear neutrophils \cite{behr2015localized}. The attraction and repulsion forces yield oscillatory trajectories for both membranes, shown in Fig \ref{fig: tctrajsnaps}. Table \ref{tab: tcspeedup} shows the increase in the speedup factor as the space grid becomes finer.

\begin{figure}[ht]
    \centering
    \includegraphics[width=1\textwidth]{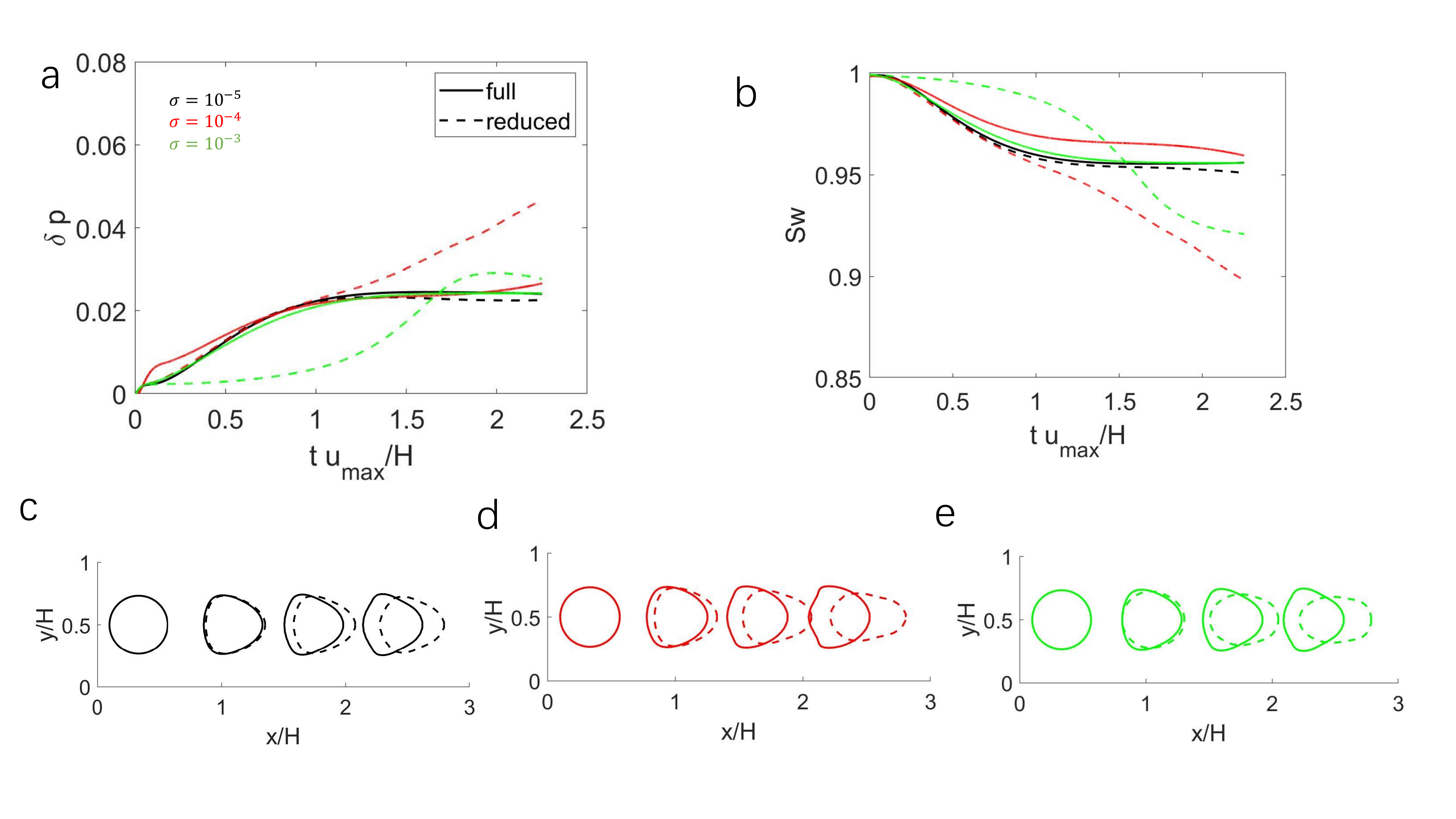}
    \caption{Transport of a circular capsule in plane-Poiseuille flow. (a) Variation of the capsule relative perimeter over time as function of the spring constant $\sigma$. (b) Variation of the capsule swelling ratio over time as function of the spring constant $\sigma$. (c,d,e) Capsule snapshots for $\sigma = 10^{-5}$ (c), $\sigma = 10^{-4}$ (d), and $\sigma = 10^{-3}$ (e). (Solid curves for FOM and dashed curves for ROM)}
    \label{fig: capsule_trans}
\end{figure}

{\color{black}
\subsection{Transport of circular capsule in a plain-Poiseuille flow}
In this test case, the dynamics of a capsule within a plane-Poiseuille flow is considered. The setup of this example follows the test conducted by Coclite et al.\cite{COCLITE201941} Initially, the capsule has a diameter of $7$ $\mu m$ and is immersed in a $2D$ channel with a height $H = 15  \ \mu m$ and length equal $3H$, centered at $7.5$ $\mu m$ away from the bottom of the lower wall. The fluid is initially at rest, the plane-Poiseuille flow with $u_{max}=10 \ \mu m/s$ is then established by posing a linear pressure drop. Simulation is run at $Re=0.01$, with $\rho = 100  \ kg\cdot m^{-3}$, $\mu = 10^-5 Pa\cdot s$. The body force on the capsule is the same as in Section 4.1, with three spring constant $\sigma=10^{-5} N\cdot \mu m^{-1}$, $10^{-4} N \cdot \mu m^{-1}$ and $10^{-3} N\cdot \mu m^{-1}$.

Following Coclite et al. \cite{COCLITE201941}, we compare the results between the FOM and the ROM in terms of the capsule perimeter variation with respect to its original configuration, $\delta p(t) = \frac{p(t)-p_0}{p_0}$ (Fig \ref{fig: capsule_trans}a), and of the swelling ratio, $Sw=\frac{A(t)}{p^2(t)/4\pi}$, where $A(t)$ is the area associated with a circle of perimeter $p(t)$ (Fig \ref{fig: capsule_trans}b). The snapshots of FOM and ROM are also compared (Fig \ref{fig: capsule_trans}c, d \& e).

For $\sigma = 10^{-5}$, the ROM is a fair approximation of the FOM. As the force coefficient increases, the system becomes more stiff. Consequently, the ROM simulation does not approximate the FOM well. We emphasize that our result is not in full agreement with the data published in Coclite et al. \cite{COCLITE201941} for two reasons. First, the time-dependent Stokes equations are considered in this work instead of Navier-Stokes equations. Secondly, the force we applied to the cell model is different.}

\section{Conclusion}\label{sec:con}
In this paper, we develop a reduced-order modeling framework for  FSI problems. Using the IBM as an example, we discussed the transfer function and its approximations.
This proposed ROM formulation enforces the impressibility condition and also preserves the Lyapunov stability. An efficient interpolation technique is applied to efficiently update the time-dependent coefficient matrices. 
The proposed model reduction technique
is applied to several biological applications {\color{black}involving linear incompressible Stokes flows}, as demonstrated
by the examples.
Compared to other traditional methods,
this new method has the following two advantages: 1) the fluid variables are  the most time-consuming part in the traditional methods, such as IBM, IIM, and FDM. But they are not explicitly involved in our ROM; 2) the structure equation is derived explicitly. It does not require special  discretization techniques, e.g., those for singular integrals used in the BEM. 
Recently, there have been growing interest in combining the reduced-order technique and data-driven methods. In this scenario, rather than the direct access to the FOM, one works with observations, e.g., structure conformations, in the form of time series. The problem is then reduced to inferring parameters in the ROM. This work is underway.    

\section*{Acknowledgments}
This work is supported by the National Science Foundation Grants  DMS-1953120 (XL) and DMS-2052685 (WH).



\begin{table}[htb]
    \caption{Membrane oscillation: Speedup of full order model and reduced-order model.}
    \label{tab: memspeedup}
    \centering
    \begin{tabular}{|l|r|r|r|l|c|}
    \hline
    \multicolumn{1}{|c|}{\multirow{2}{*}{h}} & \multicolumn{2}{c|}{Model order}  & \multicolumn{2}{c|}{CPU time}  & \multirow{2}{*}{Speedup factor} \\ \cline{2-5}
    \multicolumn{1}{|c|}{}    & \multicolumn{1}{c|}{full} & \multicolumn{1}{c|}{reduced} & \multicolumn{1}{c|}{full} & \multicolumn{1}{c|}{reduced} &  \\ \hline
    1/6   & 1728  & 144   & .0118  & .0036     & 3.2778   \\ \hline
    1/8   & 3072  & 192   & .0309  & .0044     & 7.0227   \\ \hline
    1/12  & 6912  & 288   & .1391  & .007      & 19.871   \\ \hline
    1/16  & 12288 & 384   & .3940  & .0166     & 24.735   \\ \hline
    1/20  & 19200 & 480   & .9745  & .0275     & 35.436   \\ \hline
    \end{tabular}
\end{table}

\begin{table}[htb]
    \caption{Membrane oscillation: Sampling cost and overall expected time saving (in seconds) for different numbers of total time steps.}
    \label{tab: memcost}
    \centering
    \begin{tabular}{|c|c|c|c|c|c|c|c|c|c|c|}
    \hline
    \multicolumn{1}{|c}{\multirow{2}{*}{h}} & \multicolumn{1}{|c}{\multirow{2}{*}{Sampling time}} &
    \multicolumn{3}{|c|}{$N_T=15$} & 
    \multicolumn{3}{|c|}{$N_T=30$} &
    \multicolumn{3}{|c|}{$N_T=50$} \\
    \cline{3-11}
    \multicolumn{1}{|c|}{}    & \multicolumn{1}{c}{} & \multicolumn{1}{|c}{FOM} & \multicolumn{1}{|c|}{ROM} & \multicolumn{1}{c|}{Saving}
    & \multicolumn{1}{|c}{FOM} & \multicolumn{1}{|c|}{ROM} & \multicolumn{1}{c|}{Saving}
    & \multicolumn{1}{|c}{FOM} & \multicolumn{1}{|c|}{ROM} & \multicolumn{1}{c|}{Saving}\\ \hline
    1/8   & .053  & .464 & .119 & .345 & .927 & .185 & .742 & 1.55 & .273 & 1.28 \\ \hline
    1/12  & .170  & 2.09 & .275 & 1.81 & 4.17 & .380 & 3.79 & 6.95 & .520 & 6.43 \\ \hline
    1/16  & .748 & 5.91 & .997 & 4.94 & 11.8 & 1.246 & 10.5 & 19.7 & 1.58 & 18.1  \\ \hline
    \end{tabular}
\end{table}

\begin{table}[htb]
    \caption{Particle rotation: Speedup of full order model and reduced-order model.}
    \label{tab: orbspeedup}
    \centering
    \begin{tabular}{|l|r|r|r|l|c|}
    \hline
    \multicolumn{1}{|c|}{\multirow{2}{*}{h}} & \multicolumn{2}{c|}{Model order}  & \multicolumn{2}{c|}{CPU time}  & \multirow{2}{*}{Speedup factor} \\ \cline{2-5}
    \multicolumn{1}{|c|}{}    & \multicolumn{1}{c|}{full} & \multicolumn{1}{c|}{reduced} & \multicolumn{1}{c|}{full} & \multicolumn{1}{c|}{reduced} &  \\ \hline
    3/16   & 2048  & 32   & .0214  & .0017     & 12.5882   \\ \hline
    1/8   & 4608  & 48   & .1012  & .0029     & 34.8966   \\ \hline
    3/32  & 8192  & 64  &  .3065 &  .0057     & 53.7719   \\ \hline
    \end{tabular}
\end{table}

\begin{table}[htb]
    \caption{Two cells interaction: Speedup of full order model and reduced-order model.}
    \label{tab: tcspeedup}
    \centering
    \begin{tabular}{|l|r|r|r|l|c|}
    \hline
    \multicolumn{1}{|c|}{\multirow{2}{*}{h}} & \multicolumn{2}{c|}{Model order}  & \multicolumn{2}{c|}{CPU time}  & \multirow{2}{*}{Speedup factor} \\ \cline{2-5}
    \multicolumn{1}{|c|}{}    & \multicolumn{1}{c|}{full} & \multicolumn{1}{c|}{reduced} & \multicolumn{1}{c|}{full} & \multicolumn{1}{c|}{reduced} &  \\ \hline
    3/16   & 4096  & 64   & .0834  & .0080     & 10.425   \\ \hline
    1/8   & 9216  & 96   & .4272  & .0148     & 28.8649   \\ \hline
    3/32  & 16384  & 128  &  1.2492 &  .0378     & 33.0476   \\ \hline
    \end{tabular}
\end{table}

\bibliographystyle{unsrt}  
\bibliography{reduced-order}

\end{document}